\documentclass[11pt,a4paper]{article}

\usepackage[utf8]{inputenc}
\usepackage[T1]{fontenc}
\usepackage{lmodern}
\usepackage[margin=2.6cm]{geometry}
\usepackage{amsmath,amssymb,amsthm,mathtools}
\usepackage{booktabs}
\usepackage{graphicx}
\graphicspath{{figures/}}
\usepackage{enumitem}
\usepackage[numbers,sort&compress]{natbib}
\usepackage[colorlinks=true,linkcolor=blue!50!black,citecolor=green!40!black,urlcolor=blue!50!black]{hyperref}
\usepackage{xcolor}

\theoremstyle{plain}
\newtheorem{theorem}{Theorem}[section]
\newtheorem{lemma}[theorem]{Lemma}
\newtheorem{proposition}[theorem]{Proposition}
\newtheorem{corollary}[theorem]{Corollary}

\theoremstyle{definition}
\newtheorem{definition}[theorem]{Definition}
\newtheorem{example}[theorem]{Example}
\theoremstyle{remark}
\newtheorem{remark}[theorem]{Remark}

\newcommand{\R}{\mathbb{R}}
\newcommand{\eps}{\varepsilon}
\newcommand{\ce}[1]{\mathrm{#1}}
\DeclareMathOperator{\rank}{rank}
\DeclareMathOperator{\supp}{supp}
\DeclareMathOperator{\sign}{sign}

\title{Canard points in planar mass-action networks:\\
a flux criterion, a duality, and a complete classification}
\author{Miruna-Mihaela Beldiman\thanks{Institute of Mathematical Statistics
and Applied Mathematics of the Romanian Academy, Bucharest, Romania.
ORCID 0000-0001-9947-7594.}
\and
Andrei-Dan Halanay\thanks{University of Bucharest, Bucharest, Romania.
ORCID 0009-0005-5481-8021.}}
\date{}

\begin{document}
\maketitle

\begin{abstract}
We determine which mass-action reaction networks on two species, one fast
and one slow, admit a canard point: a singular Hopf bifurcation at a fold of
the critical manifold, followed by a canard explosion. In flux coordinates,
existence reduces to three linear equations, one linear inequation and a sign
condition on the product of two linear forms, and the criticality of the
explosion is the sign of an explicit cubic form. For four reactions with a
balance matrix of full rank, both become integer polynomial conditions on the
stoichiometry. A vacancy duality, which exchanges the activator--inhibitor
and substrate-depletion mechanisms, preserves canard points together with
their criticality. The central structural result is a mode criterion: a
network has a canard point if and only if its elementary flux modes contain
a sign-meeting pair and a mode carrying a fold. Hence every canard network
contains a minimal one carried by at most three modes, and so by at most
twelve reactions, whatever the molecularity. For complexes of molecularity at
most three we enumerate all 138\,754 minimal canard networks and label each
exactly by the criticalities it admits. None has more than six reactions; 32
have canard points of codimension zero and no explosion; and most
six-reaction networks admit both criticalities, separated by degenerate
canard points. All labels rest on certificates checked in exact rational
arithmetic. Resolving the cubic steps into bimolecular ones through fast
dimers keeps the canard points and changes the criticality only by an
explicit term.
\end{abstract}

\medskip
\noindent\textbf{Keywords:} canard explosion, chemical reaction networks,
singular Hopf bifurcation, elementary flux modes.

\smallskip
\noindent\textbf{MSC 2020:} 34E17 (primary); 34C23, 37G15, 34E15, 92C42,
80A30 (secondary).

\section{Introduction}

\paragraph{Canards.}
A canard explosion is the transition from small limit cycles, born in a Hopf
bifurcation, to relaxation oscillations. The intermediate cycles, the canard
cycles, exist only within an exponentially small interval of a parameter, so
the transition is abrupt. It occurs in slow--fast systems at a
\emph{canard point}, a fold of the critical manifold at which the slow flow
passes from the attracting to the repelling branch. In planar systems it
occurs generically in one-parameter families that unfold a canard point, and
it is the way in which relaxation oscillations of van der Pol type are born
\citep{KrupaSzmolyan2001relax}.

The phenomenon was discovered with nonstandard analysis
\citep{BenoitCallotDiener1981}. It was later explained geometrically by
\citet{KrupaSzmolyan2001fold,KrupaSzmolyan2001relax}, who used the blow-up
method. Their work provides the two ingredients we rely on: the asymptotic
location of the explosion, and a coefficient $A$ whose sign decides whether
the singular Hopf bifurcation, and with it the explosion, is supercritical or
subcritical. In three dimensions, folded singularities and their canards
\citep{SzmolyanWechselberger2001,Wechselberger2005} organise mixed-mode
oscillations \citep{Desroches2012}. The theory is surveyed in the monographs
of \citet{Kuehn2015}, \citet{Wechselberger2020} and
\citet{DeMaesschalckDumortierRoussarie2021}, and the blow-up method in
\citet{JardonKojakhmetovKuehn2021}. \citet{DeMaesschalckDoanWynen2021} gave an
intrinsic, coordinate-free criterion for the criticality.

\paragraph{Canards in chemical kinetics.}
Chemical reactions were among the first sources of canards and mixed-mode
oscillations. Canard explosions were found in the Oregonator model of the
Belousov--Zhabotinsky reaction \citep{BronsBarEli1991,BronsBarEli1994}; its
two-variable reduction belongs to the class classified here
(Section~\ref{sec:disc}). Mixed-mode oscillations were studied in chemical systems by
\citet{PetrovScottShowalter1992} and \citet{MilikSzmolyan2001}. Autocatalator
models show canard explosions and relaxation oscillations
\citep{BronsSturis2001,GucwaSzmolyan2009} and, in three dimensions, twin
canard orbits and mixed-mode oscillations \citep{HasanKrauskopfOsinga2017}.
Further examples are the Olsen model of the peroxidase--oxidase reaction
\citep{KuehnSzmolyan2015}, the Brusselator \citep{EngelOliconMendez2026}, the
Sel'kov oscillator \citep{BrechmannRendall2018} and, beyond chemistry, a
predator--prey model with an Allee effect \citep{LiLiWang2023}. In each case
the analysis is carried out for one model, with its own coordinates and
scaling. The question we ask is a structural one: \emph{which reaction
networks admit a canard explosion, and of which criticality?}

\paragraph{Time scales in reaction networks.}
To ask the question for networks, one needs slow--fast splittings that come
from the network. Two approaches identify small parameters and the resulting
reductions systematically. The first uses Tikhonov--Fenichel parameter values
\citep{GoekeWalcher2014,GoekeWalcherZerz2015,FeliuKruffWalcher2020,FeliuWalcherWiuf2022},
which \citet{Kruff2021} compute algorithmically. The second uses tropical
scalings of rate constants and concentrations
\citep{RadulescuVakulenkoGrigoriev2015,DesoeuvresSzmolyanRadulescu2022}; for
planar polynomial systems see \citet{KristiansenSarantaris2026}. We take the
splitting as given: every change of one species is slow. This is the form
such scalings produce when that species varies on the slow time scale.

\paragraph{Structural bifurcation theory.}
Chemical reaction network theory relates dynamics to network structure
\citep{Feinberg2019}. Small networks have been classified with respect to
multistationarity
\citep{JoshiShiu2017,ShiuDeWolff2019,TangLinZhang2022,TangZhang2022}, and
parameter regions for multistationarity have been described by the sign of a
Jacobian determinant on the steady-state variety
\citep{ConradiFeliuMinchevaWiuf2017}. The smallest bimolecular networks
with a Hopf bifurcation were determined by \citet{BanajiBoros2023}, and the
bifurcations and oscillations of small planar networks were studied by
\citet{BanajiBorosHofbauer2024b} and \citet{BorosHofbauer2024}. There are
structural conditions for saddle-node bifurcations \citep{Vassena2023}, and
folds and cusps of mass-action networks can be computed
\citep{Banaji2026fold}. Finally, oscillations \citep{Banaji2018} and, more
generally, nondegenerate local bifurcations \citep{BanajiBorosHofbauer2025}
are inherited when a network is enlarged; see \citet{BanajiBorosHofbauer2025}
and the references therein for other inheritance results.

A canard point combines three requirements: a nondegenerate fold of the fast
nullcline, a balance of the slow flux at that fold, and a sign condition that
makes the bifurcation a Hopf bifurcation rather than a saddle. Its
criticality is a cubic quantity. To our knowledge, neither has been treated
structurally. Two features of our approach should be noted. First, several
of the existing criteria are proved for \emph{parameter-rich} kinetics, in
which the derivatives of the rates can be chosen independently, and, as
\citet{Vassena2023} and \citet{VassenaStadler2024} stress, mass action is not
of this kind. We work directly with mass action, in the convex coordinates of
\citet{Clarke1980}. Second, we use elementary flux modes, classical in the
analysis of metabolic and reaction networks \citep{SchusterHilgetag1994}, to
reduce the canard question to configurations of at most three modes.

\paragraph{Setting.}
Let X and Y be the species. Reaction $r$ has reactant
$a_r\ce X+b_r\ce Y$ and changes X and Y by $N_r$ and $M_r$. We consider the
slow--fast mass-action system
\begin{equation}\label{eq:sys}
  \dot x=f(x,y)=\sum_rN_rk_rx^{a_r}y^{b_r},\qquad
  \dot y=\eps\beta(x,y)=\eps\sum_rM_rk_rx^{a_r}y^{b_r},
\end{equation}
in which every change of Y is slow. Systems of this form arise from mass-action
networks after a scaling of rate constants and concentrations under which
the Y-changes become slow \citep{GoekeWalcherZerz2015,FeliuKruffWalcher2020}.

\paragraph{Main results.}
The first result is a \emph{flux criterion} (Theorem~\ref{thm:planar}).
Write $v_r=k_rx^{a_r}y^{b_r}$ for the fluxes, and define the moment forms
\[
  n_{ij}=\textstyle\sum_rN_r(a_r)_i(b_r)_jv_r,\qquad
  m_{ij}=\sum_rM_r(a_r)_i(b_r)_jv_r,
\]
where $(a)_i$ is the falling factorial. A canard point exists for some rate
constants if and only if some $v>0$ satisfies
\[
  n_{00}=n_{10}=m_{00}=0,\qquad n_{20}\neq0,\qquad n_{01}m_{10}<0 ,
\]
and the explosion is supercritical if and only if the cubic $Q$
of~\eqref{eq:Q} is negative. For four reactions whose balance matrix has full
rank, the flux is determined by $3\times3$ minors, so existence and
criticality become integer polynomial sign conditions on the stoichiometry
(Corollary~\ref{cor:four}). The map
$a\ce X+b\ce Y\mapsto a\ce X+(d-a-b)\ce Y$ is a \emph{vacancy duality}
(Theorem~\ref{thm:dual}): it preserves canard points, fluxes and
criticality, and exchanges activator--inhibitor with substrate depletion. It
comes from the change of coordinates $(x,y)\mapsto(x/y,1/y)$, a change of
chart once complexes are padded with vacancies
(Proposition~\ref{prop:chart}).

The central result is the \emph{mode criterion} (Theorem~\ref{thm:modes}): a
network has a canard point if and only if its elementary flux modes contain
a sign-meeting pair and a mode carrying a fold. Elementary modes are the
minimal steady-state pathways of the network, as in metabolic pathway
analysis \citep{SchusterHilgetag1994}. A mode carries a fold if the fast
nullcline bends along it; such a mode contains a cooperative step. Two modes
are sign-meeting if, mixed together, they supply the slow feedback with one
of the two signs that give a Hopf bifurcation rather than a saddle. Minimal
canard networks are therefore carried by at most three modes, i.e.\ by at
most twelve reactions (Corollary~\ref{cor:three}). This reduces their enumeration to the elementary
modes. In Section~\ref{sec:class} we give the complete list for molecularity
at most three, with exact criticality labels (Table~\ref{tab:class}). All
three mode types occur, and no minimal network has more than six reactions.

We also show that a fold of the fast nullcline requires an X-producing
reaction with at least two X in its reactant (Corollary~\ref{cor:coop}). For
molecularity at most three, this means $2\ce X\to3\ce X$ or
$2\ce X+\ce Y\to3\ce X$. The underlying fold criterion
(Proposition~\ref{prop:fold}) is equivalent to the 2-alternating condition of
\citet{JoshiShiu2017}. We present it as the fold counterpart of their
multistationarity classification, not as a new classification. Steps with
three molecules can be resolved into bimolecular ones through fast dimers.
This keeps the canard points, and the criticality is then the sign of
$Q+2\theta n_{01}m_{10}n_{20}$, where $\theta\in(0,1)$ measures the
strength of the dimerisation (Corollary~\ref{cor:bimol}).

\paragraph{What is new.}
The flux criterion with its explicit cubic form for the criticality, the
vacancy duality, the mode criterion, the classification and the correction
of the criticality under fast storage are new. The mode
criterion bounds minimal canard networks by twelve reactions independently of
the molecularity. The classification shows that minimal canard networks need
not be small, since here they have up to six reactions, and that many of them
admit both criticalities. It is computational but certified: every label
rests on an exact certificate, a finite object checked in rational
arithmetic. A one-sided label is certified by coefficients of $Q$ of one sign
over the extreme rays of the admissible cones, and the label \emph{both} by
two rational fluxes at which $Q$ takes opposite signs. Minimality was also
checked on random samples by linear programming, independently of the
elementary modes, and every entry was re-derived by an independent program in
exact arithmetic (Section~\ref{sec:class}). The list, the certificates, the
code and the logs of all computations are openly available
\citep{CanardPlanarData2026}.

\paragraph{Relation to geometric singular perturbation theory.}
The dynamical input is the theorem of \citet{KrupaSzmolyan2001relax}, and
nothing in our analysis replaces it. What we add is that, for a network, the
set of rate constants and points at which its hypotheses hold becomes a
polyhedral object. In flux coordinates the fold, the balance and the
unfolding conditions are linear, and the criticality is the sign of a cubic
form. The studies of chemical canards cited above verify the hypotheses for
one model, in coordinates adapted to it, and the criterion of
\citet{DeMaesschalckDoanWynen2021} gives the criticality of a given system
intrinsically. Our criterion decides existence and criticality for all rate
constants of a network at once, and this is what makes a classification
possible.

\paragraph{Organisation.}
Section~\ref{sec:setting} recalls canard points and the coefficient $A$.
Section~\ref{sec:flux} proves the flux criterion, and Section~\ref{sec:fold}
treats the fold. Section~\ref{sec:four} specialises to four reactions, and
Section~\ref{sec:dual} introduces the vacancy duality.
Section~\ref{sec:modes} proves the mode criterion, Section~\ref{sec:class}
presents the classification and its checks, and Section~\ref{sec:bimol}
realises the networks by bimolecular steps. Section~\ref{sec:disc} distils
design rules, highlights robust networks, discusses molecularity and
thermodynamics, places the Oregonator in the classification, and lists open
questions.

\section{Canard points and their criticality}\label{sec:setting}

The critical manifold of~\eqref{eq:sys} is $\{f=0\}$. It is attracting
where $f_x<0$ and repelling where $f_x>0$.

\begin{definition}\label{def:canard}
A point $p=(x,y)\in\R^2_{>0}$ is a \emph{(nondegenerate) canard point} of
Hopf type if all of the following hold at $p$:
\begin{itemize}[nosep]
  \item $f=f_x=\beta=0$;
  \item $f_{xx}\neq0$, $f_y\neq0$ and $\beta_x\neq0$;
  \item $f_y\beta_x<0$.
\end{itemize}
\end{definition}

These are the conditions of \citet{KrupaSzmolyan2001relax} at a canard
point, after reflections of $x$ and $y$. The sign condition $f_y\beta_x<0$
makes the equilibrium near $p$ a focus or node rather than a saddle. Let
$\lambda$ be a parameter that unfolds $\beta$ at the fold. Then, for small
$\eps>0$, there is a singular Hopf bifurcation at
$\lambda_H(\sqrt\eps)$ followed by a canard explosion at
$\lambda_c(\sqrt\eps)$. Its criticality is given by the sign of the
coefficient $A$ of \citet{KrupaSzmolyan2001relax}: the explosion is
supercritical if $A<0$ and subcritical if $A>0$.

\begin{lemma}\label{lem:A}
At a canard point, $A=2\mathfrak B/\sqrt{-f_y\beta_x}$, where
\[
  \mathfrak B\,f_{xx}^2=(\beta_xf_{xy}-\beta_yf_{xx})f_{xx}+f_y(\beta_{xx}f_{xx}-\beta_xf_{xxx}).
\]
\end{lemma}

\begin{proof}
Translate $p$ to the origin and rescale $x=\alpha X$, $y=bY$, $t=\tau T$,
with
\[
  b=\frac{2\beta_x}{f_{xx}},\qquad \alpha^2=-\frac{4f_y\beta_x}{f_{xx}^2},\qquad
  \tau=\frac{2}{\alpha f_{xx}},
\]
choosing the sign of $\alpha$ so that $\tau>0$. Then
$X'=-Y+X^2+\cdots$ and $Y'=\eps(X+\cdots)$. Reading off the coefficients
$a_2,\dots,a_5$ of the canonical form of \citet{KrupaSzmolyan2001relax} and
substituting into $A=-a_2+3a_3-2a_4-2a_5$ gives the formula. The remaining
normalisations, of $\lambda$ and an $O(\eps)$ shift of $Y$, do not change
$a_2,\dots,a_5$. The same formula follows from the intrinsic criterion of
\citet{DeMaesschalckDoanWynen2021}.
\end{proof}

The coefficients $a_2,\dots,a_5$ are evaluated at $\eps=0$, and
\citet[(3.1)]{KrupaSzmolyan2001relax} allow the vector field to depend on
$\eps$. So Definition~\ref{def:canard} and Lemma~\ref{lem:A} apply
unchanged to a fast equation $\dot x=f+\eps f_1$ with $f_1$ smooth: the term
$\eps f_1$ enters only $h_3$ of the canonical form, after an $O(\eps)$ shift
of $y$ if $f_1(p)\neq0$.

\begin{remark}[Location of the explosion]\label{rem:S}
Let $\mu$ be an unfolding parameter, and let
\[
  G_\mu=\tfrac{d}{d\mu}\beta\big(p_F(\mu);\mu\big)
\]
be the derivative of $\beta$ along the continued fold point $p_F(\mu)$.
Suppose $G_\mu\neq0$.
\citet{KrupaSzmolyan2001relax} show that
$\lambda_c-\lambda_H=-\tfrac18A\eps+O(\eps^{3/2})$. The rescaling in the proof
of Lemma~\ref{lem:A} gives $\lambda=-G_\mu(\mu-\mu_0)/(\beta_x\alpha)$, up to
higher-order terms and an $O(\eps)$ shift that is common to $\lambda_H$ and
$\lambda_c$. Hence
\[
  \mu_c(\eps)-\mu_H(\eps)=S\eps+O(\eps^{3/2}),\qquad
  S=\frac{\mathfrak B\,\beta_x}{2f_{xx}\,G_\mu}.
\]
The explosion lies on the side of $\mu_H$ where the Hopf cycles exist.
\end{remark}

\begin{remark}[What is observed for finite $\eps$]\label{rem:finite}
The exponentially small interval of the explosion contains the canard
cycles, which interpolate between the small Hopf cycles and the relaxation
oscillation. By Remark~\ref{rem:S}, it lies at distance $S\eps+O(\eps^{3/2})$
from the singular Hopf point. So it is the transition that is sharp; there is
no exponentially small window that a parameter sweep has to hit. Suppose,
as in the global setting of \citet{KrupaSzmolyan2001relax}, that the
critical manifold is S-shaped, so that relaxation oscillations exist beyond
the explosion. For fixed small $\eps$, a slow sweep of $\mu$ then shows one
of two scenarios, and the sign of $A$ decides which. If $A<0$, stable small cycles grow from $\mu_H$ over an
interval of length about $|S|\eps$, and their amplitude then jumps to that of
the relaxation oscillation: the onset of oscillations is soft. If $A>0$, the
cycles born at $\mu_H$ are unstable, and the stable equilibrium coexists with
a stable relaxation cycle over an interval of length about $|S|\eps$: the
onset is hard, and oscillations start and stop at different parameter values.
Figure~\ref{fig:explosion} shows both scenarios at $\eps=0.01$. The
criticality labels of Section~\ref{sec:class} are therefore predictions of
observable behaviour.
\end{remark}

\begin{remark}[How small $\eps$ must be]\label{rem:BT}
At a Hopf point $f_x=-\eps\beta_y$, so the Jacobian has determinant
$-\eps(\eps\beta_y^2+f_y\beta_x)$, and a Hopf point needs
$\eps\beta_y^2<-f_y\beta_x$ there. As $\eps$ grows, this inequality can fail
along the Hopf branch that starts at the canard point. The branch then ends
at a point with a double zero eigenvalue, generically a Bogdanov--Takens
point. It can also end at a fold of the Hopf curve, or persist. Evaluated at
the canard point, the inequality gives the scale
\[
  \eps_{\mathrm{BT}}=-\frac{f_y\beta_x}{\beta_y^2}=-\frac{y}{x}\,\frac{n_{01}m_{10}}{m_{01}^2}.
\]
This is an estimate, not a bound. The Hopf point lies at distance $O(\eps)$
from the canard point, which changes $f_y\beta_x$ by a term of order $\eps$
that depends on the unfolding parameter. The estimate is accurate when this
term is small compared with $\eps\beta_y^2$, that is, when the
self-regulation of the slow species, measured by $m_{01}=y\beta_y$, is
strong compared with the feedback loop $-n_{01}m_{10}$, and $\eps_{\mathrm{BT}}$
is small. Then the canard scenario needs roughly
$\eps\,m_{01}^2<-(y/x)\,n_{01}m_{10}$, and Section~\ref{sec:class} confirms
this numerically. Otherwise the estimate can fail. For the autocatalator of
Example~\ref{ex:auto} at $k=(2,1,2,1)$ we have $\eps_{\mathrm{BT}}=1$, but the
determinant is $\eps k_3(k_2+k_4x^2)>0$ everywhere, and, unfolded by $k_2$,
the Hopf point exists for every $\eps<2$.
\end{remark}

\section{The flux criterion}\label{sec:flux}

At a point $(x,y)>0$, the fluxes $v_r=k_rx^{a_r}y^{b_r}$ can take any
positive value, since $k_r=v_r/(x^{a_r}y^{b_r})$. With
$(a)_i=a(a-1)\cdots(a-i+1)$, define the moment forms
\[
  n_{ij}(v)=\sum_rN_r(a_r)_i(b_r)_j\,v_r,\qquad m_{ij}(v)=\sum_rM_r(a_r)_i(b_r)_j\,v_r .
\]
Since $x^i\partial_x^i\,x^a=(a)_ix^a$, we have
\begin{equation}\label{eq:moments}
  x^iy^j\,\partial_x^i\partial_y^jf=n_{ij}(v),\qquad x^iy^j\,\partial_x^i\partial_y^j\beta=m_{ij}(v).
\end{equation}
Let $B=[N;\,N\circ a;\,M]$ be the $3\times R$ \emph{balance matrix}, and
$K=\{v\in\R^R_{>0}:Bv=0\}$ the open \emph{canard cone}. Define
\begin{equation}\label{eq:Q}
  Q=(m_{10}n_{11}-m_{01}n_{20})\,n_{20}+n_{01}\,(m_{20}n_{20}-m_{10}n_{30}).
\end{equation}

\begin{theorem}[Flux criterion]\label{thm:planar}
\leavevmode
\begin{enumerate}[label=(\alph*),nosep]
  \item For given rate constants, $p$ is a canard point iff its fluxes $v$
        satisfy
        \begin{equation}\label{eq:cone}
          v\in K,\qquad n_{20}(v)\neq0,\qquad n_{01}(v)\,m_{10}(v)<0 .
        \end{equation}
  \item Hence some rate constants give a canard point iff~\eqref{eq:cone}
        has a solution.
  \item At a canard point, $\sign A=\sign Q(v)$.
  \item Suppose $\rank B=3$. Near each canard point, the rate constants with
        a canard point form a smooth hypersurface, and some single rate
        constant is transverse to it, i.e.\ unfolds the canard point.
        Suppose instead $\rank B<3$. Then the set of such rate constants has
        nonempty interior (codimension zero).
\end{enumerate}
\end{theorem}

\begin{proof}
(a) By~\eqref{eq:moments}:
\begin{itemize}[nosep]
  \item $f=n_{00}$, $xf_x=n_{10}$ and $\beta=m_{00}$;
  \item $x^2f_{xx}=n_{20}$;
  \item $yf_y=n_{01}$ and $x\beta_x=m_{10}$.
\end{itemize}
Positive factors do not change signs.

(b) This follows from~(a), since every $v>0$ is realised at every point.

(c) By~\eqref{eq:moments}, $x^4y\,\mathfrak B\,f_{xx}^2=Q(v)$, and
$\sign A=\sign\mathfrak B$ by Lemma~\ref{lem:A}.

(d) Fix $k$. The Jacobian of $(f,f_x)$ with respect to $(x,y)$ has
determinant $f_xf_{xy}-f_yf_{xx}=-f_yf_{xx}\neq0$ at a canard point. So the
fold point $p_F(k)$ is locally a smooth function of $k$, and the canard set is
$\{g(k)=\beta(p_F(k),k)=0\}$. The derivative of $(f,f_x,\beta)$ with respect
to $k$ at fixed $(x,y)$ is
\[
  \operatorname{diag}(1,1/x,1)\,B\,\operatorname{diag}(x^{a_r}y^{b_r}).
\]

If $\rank B=3$, choose $\delta k$ with $D_k(f,f_x)\,\delta k=0$ and
$D_k\beta\,\delta k=1$. Then $p_F$ does not move to first order and
$Dg\cdot\delta k=1$. So $\nabla g\neq0$, and some $\partial g/\partial k_r\neq0$.

If $\rank B<3$, the solutions of~\eqref{eq:cone} in the coordinates
$(x,y,v)$ form a manifold of dimension $2+R-\rank B\ge R$. Its map to rate
constants $k_r=v_r/(x^{a_r}y^{b_r})$ has discrete fibres, because for fixed
$k$ the canard points are isolated: the Jacobian above is invertible. By the
constant-rank theorem, on the open set where the rank of this map is maximal
the rank must be $R$, so the image has nonempty interior.
\end{proof}

We call the canard points of networks with $\rank B<3$
\emph{codimension-zero}. For a network with a canard point, the rows $N$
and $N\circ a$ are independent. Otherwise $N\circ a=cN$, so
$n_{20}=c(c-1)n_{00}=0$ on $K$. Hence $\rank B\ge2$. If $\rank B=2$, then $M$
is a combination of $N$ and $N\circ a$, so $\beta=\alpha_1 f+\alpha_2 xf_x$,
and the fold point is an equilibrium for every choice of rate constants.

Its Jacobian there has determinant $-\eps f_y\beta_x>0$ and trace
$\eps\beta_y=\eps m_{01}/y$. For the $32$ minimal codimension-zero networks
of Section~\ref{sec:class}, exact linear programming shows that $m_{01}\le0$
on the admissible set:
\begin{itemize}[nosep]
  \item for $27$ of them the trace is negative except on a lower-dimensional
        set, so the equilibrium is stable and no Hopf bifurcation occurs;
  \item for the other $5$, $m_{01}\equiv0$, and the equilibrium is a linear
        centre for all rate constants.
\end{itemize}
In neither case does a canard explosion occur within these networks. They
matter as building blocks.

The sign condition admits two mechanisms:
\begin{itemize}[nosep]
  \item \emph{activator--inhibitor} (AI): $f_y<0<\beta_x$. Y inhibits the
        growth of X, and X raises Y.
  \item \emph{substrate depletion} (SD): $f_y>0>\beta_x$. Y feeds X, and X
        consumes Y. The autocatalator and the Brusselator are of this kind.
\end{itemize}

\section{The fold requires cooperativity}\label{sec:fold}

The first two equations of~\eqref{eq:cone}, $n_{00}=n_{10}=0$, together with
$n_{20}\neq0$, say that the fast nullcline has a nondegenerate fold. Let
$P=\{a_r:N_r>0\}$ and $C=\{a_r:N_r<0\}$ be the X-orders of the producing and
the consuming reactions. Note that $C\subset\{1,2,\dots\}$. If $P$ or $C$ is
empty, then $n_{00}=0$ has no positive solution. Only reactions with
$N_r\neq0$ matter here, and the problem is that of the one-species network
obtained by projecting onto X.

\begin{proposition}\label{prop:fold}
There are $v>0$ with $n_{00}=n_{10}=0\neq n_{20}$ iff the following hold:
\begin{itemize}[nosep]
  \item $\max P>\min C$;
  \item $\min P<\max C$;
  \item $P$ and $C$ are not one and the same two-element set.
\end{itemize}
Equivalently, there are orders $a<b<c$ that are realised by producing,
consuming and producing reactions respectively, or by consuming, producing
and consuming ones. This is a $2$-alternating subnetwork in the sense of
\citet[Def.~3.4]{JoshiShiu2017}.
\end{proposition}

\begin{proof}
Set $w_r=|N_r|v_r$ and normalise the total production to $1$. Then
$n_{00}=n_{10}=0$ says that two fully supported probability vectors, on the
producing and on the consuming reactions, have the same mean order. This is
possible iff the first two conditions hold, or $P=C$ is a single point.

Since $n_{10}=0$, $n_{20}=\sum_rN_ra_r(a_r-1)v_r=\sum_rN_ra_r^2v_r$ is the
difference of the second moments.
\begin{itemize}[nosep]
  \item If one side has at least three distinct orders, its second moment
        can be varied at a fixed interior mean.
  \item If one side is a single point and the other is not, Jensen's
        inequality is strict.
  \item If both sides have two orders, the second moments are affine
        functions of the common mean. They coincide identically iff the two
        order sets coincide; otherwise they agree for at most one mean, while
        the admissible means form an open interval.
\end{itemize}
The equivalence with the alternating triple is elementary: take $c_1<p_1$
and $p_2<c_2$; either a triple can be read off directly, or
$c_1=p_2$ and $p_1=c_2$, and then any third order completes a triple.
\end{proof}

\begin{corollary}\label{cor:coop}
A fold requires an X-producing reaction with at least two X in its reactant,
hence at least three X in its product. So no network of molecularity at most
two has a canard point. For molecularity at most three, every network with a
canard point contains $2\ce X\to3\ce X$ or $2\ce X+\ce Y\to3\ce X$.
\end{corollary}

\begin{proof}
$\max P>\min C\ge1$ forces a producing order of at least $2$.
\end{proof}

\citet[Thm.~3.6]{JoshiShiu2017} prove that a one-species network is
nondegenerately multistationary iff it has a 2-alternating subnetwork.
Proposition~\ref{prop:fold} is the corresponding statement for folds. The
Sel'kov model \citep{BrechmannRendall2018}, whose fast reactions are
$\ce A+2\ce B\to3\ce B$ and $\ce B\to0$ (with B fast and A slow), has
$P=\{2\}$ and $C=\{1\}$, so it has no fold. Adding $\ce A\to\ce B$, as in the autocatalator
\citep{BronsSturis2001,GucwaSzmolyan2009}, adds the producing order $0$ and
creates the fold.

\begin{remark}[Orientation of the fold]\label{rem:orient}
Let all complexes have molecularity at most three, and let
$S_k=\sum_{a_r=k}N_rv_r$ be the net production of X by the reactions with $k$
molecules of X in their reactant. Then $S_0\ge0$ is the supply of X by
reactions without X in their reactant, which cannot consume X. Every reaction
with reactant $3\ce X$ has a product with fewer X, so $N_r\le-1$, and
$S_3\le0$ is the removal of X by such reactions; $S_3<0$ exactly when one of
them is present. At every point with
$n_{00}=n_{10}=0$,
\[
  x^2f_{xx}=n_{20}=2\,(S_0+S_3).
\]
Indeed, $n_{00}=S_0+S_1+S_2+S_3=0$ and $n_{10}=S_1+2S_2+3S_3=0$ give
$S_2=S_0-2S_3$, and $n_{20}=2S_2+6S_3$. Hence a fold needs, besides the
cooperative step of Corollary~\ref{cor:coop}, a supply of X from an X-free
reactant or a removal of X at third order. At the fold, $f$ is convex in $x$
($f_{xx}>0$) when the supply outweighs the third-order removal, and it can be
concave only if a reaction with reactant $3\ce X$ is present.
\end{remark}

\section{Four reactions}\label{sec:four}

\begin{corollary}\label{cor:four}
Let $R=4$ and $\rank B=3$. Let $\kappa_j=(-1)^j\det B_{\hat j}$, where
$B_{\hat j}$ is $B$ without column $j$. Then the network has a canard point
iff three conditions hold:
\begin{itemize}[nosep]
  \item all $\kappa_j$ are nonzero and of one sign (take $\kappa>0$);
  \item $n_{20}(\kappa)\neq0$;
  \item $n_{01}(\kappa)m_{10}(\kappa)<0$.
\end{itemize}
Its criticality is $\sign Q(\kappa)$, the same for all rate constants.
\end{corollary}

\begin{proof}
$\ker B$ is spanned by $\kappa$, and $Q$ is homogeneous of degree~$3$.
\end{proof}

\begin{example}\label{ex:auto}
For the planar autocatalator
$\{0\to\ce Y,\ \ce Y\to\ce X,\ \ce X\to0,\ 2\ce X+\ce Y\to3\ce X\}$
(X fast, Y slow) we get $\kappa\propto(2,1,2,1)$. This gives
$n_{01}=2\kappa_4>0>m_{10}=-2\kappa_4$ (substrate depletion) and $Q<0$
(supercritical). This network is the autocatalator model~(1.1) of
\citet{GucwaSzmolyan2009}, written with $\ce X=b$ and $\ce Y=a$. They
established its relaxation oscillations for $\mu>1$, and its canard point
lies at $\mu=1$.
\end{example}

Among the networks of four reactions whose complexes have molecularity at
most three, exactly $4\,982$ have a canard point of codimension one
(Table~\ref{tab:four}).

\begin{table}[ht]
\centering
\small
\begin{tabular}{@{}lrrr@{}}
\toprule
& supercritical & subcritical & $A\equiv0$\\
\midrule
activator--inhibitor & $1\,230$ & $1\,157$ & $104$\\
substrate depletion & $1\,230$ & $1\,157$ & $104$\\
\bottomrule
\end{tabular}
\caption{The $4\,982$ four-reaction networks with a canard point of
codimension one ($\rank B=3$; molecularity at most three), by mechanism and
criticality.}
\label{tab:four}
\end{table}

\section{The vacancy duality}\label{sec:dual}

The symmetry of Table~\ref{tab:four} has a structural reason.

\begin{theorem}\label{thm:dual}
Let all complexes have molecularity at most $d$. The map
\[
  a\ce X+b\ce Y\ \longmapsto\ a\ce X+(d-a-b)\ce Y,
\]
applied to reactants and products, is an involution on networks. It maps a
network with a canard point at fluxes $v$ to a network with a canard point at
the same fluxes. It preserves the sign of $A$, and it exchanges
activator--inhibitor with substrate depletion.
\end{theorem}

\begin{proof}
Under the map, $a$ and $N$ are unchanged, $b\mapsto d-a-b$, and
$M\mapsto-(N+M)$. On $K$ we have $n_{00}=n_{10}=m_{00}=0$. Hence:
\begin{itemize}[nosep]
  \item $m_{00}\mapsto-(n_{00}+m_{00})$ and $m_{10}\mapsto-(n_{10}+m_{10})$,
        so $K$ is unchanged and $m_{10}\mapsto-m_{10}$;
  \item $n_{01}\mapsto dn_{00}-n_{10}-n_{01}=-n_{01}$;
  \item $n_{20}$ and $n_{30}$ are unchanged;
  \item $n_{11}\mapsto-n_{11}-n_{20}$, $m_{01}\mapsto m_{01}+m_{10}+n_{01}$ and
        $m_{20}\mapsto-(m_{20}+n_{20})$.
\end{itemize}
Substituting into~\eqref{eq:Q}, the terms $\pm n_{01}n_{20}^2$ cancel and $Q$
is unchanged. The product $n_{01}m_{10}$ is unchanged, while the sign of
$n_{01}$ flips.
\end{proof}

The duality pads every complex with ``vacancies'' up to molecularity $d$, and
exchanges Y with the vacancies. For $d=3$, the dual of the planar
autocatalator is the activator--inhibitor network
$\{3\ce Y\to2\ce Y,\ 2\ce Y\to\ce X+2\ce Y,\ \ce X+2\ce Y\to3\ce Y,\ 2\ce X\to3\ce X\}$,
which is also supercritical.

The duality is not an artefact of flux coordinates: it comes from a change
of coordinates in the state space. Write $f^*$ and $\beta^*$ for the
functions of the dual network $G^*$, with the same rate constants $k_r$, in
coordinates $(u,w)$.

\begin{proposition}[The duality as a change of chart]\label{prop:chart}
Let $\Phi(x,y)=(x/y,\,1/y)$, a diffeomorphism of $\R^2_{>0}$. After the
change of time $dt=w^{d-1}d\tau$, that is, multiplication of the vector field
by $w^{d-1}$, $\Phi$ takes the system~\eqref{eq:sys} of a network $G$ to
\begin{equation}\label{eq:chart}
  u'=f^*+\eps\,u\,(\beta^*+f^*),\qquad w'=\eps\,w\,(\beta^*+f^*).
\end{equation}
This system and the system~\eqref{eq:sys} of $G^*$ have the same canard
points, with the same sign of $A$.
\end{proposition}

\begin{proof}
Since $x^ay^b=w^{-d}u^aw^{d-a-b}$, we have
\begin{gather*}
  f(x,y)=w^{-d}f^*(u,w),\qquad \beta(x,y)=w^{-d}\hat\beta(u,w),\\
  \hat\beta=\sum_rM_rk_ru^{a_r}w^{d-a_r-b_r}=-(\beta^*+f^*),
\end{gather*}
where the last equality holds because the dual has Y-changes
$-(N_r+M_r)$. Hence
\[
  \dot u=\frac{\dot x}{y}-\frac{x\dot y}{y^2}=w^{1-d}\bigl(f^*+\eps u(\beta^*+f^*)\bigr),\qquad
  \dot w=-\frac{\dot y}{y^2}=\eps\,w^{2-d}(\beta^*+f^*),
\]
which gives~\eqref{eq:chart} after multiplication by $w^{d-1}$.

The system~\eqref{eq:chart} differs from that of $G^*$ in two ways. First,
its fast equation has an extra term of order $\eps$, which vanishes at a
canard point. By the remark after Lemma~\ref{lem:A}, it changes neither the
canard conditions nor $A$. In the canonical form of
\citet[(3.3)]{KrupaSzmolyan2001relax} it enters only $h_3$, and hence the
coefficient $a_1=\partial_xh_3(0)$, which moves the Hopf point but does not
enter $A$ \citep[(3.13)]{KrupaSzmolyan2001relax}. Second, its slow function
is $s\,(\beta^*+cf^*)$ with $s=w>0$ and $c=1$; here $s$ depends only on the
slow variable. At a point with $f^*=f^*_u=\beta^*=0$, the function
$\tilde\beta=s(\beta^*+cf^*)$ vanishes, and $\tilde\beta_u=s\beta^*_u$,
$\tilde\beta_w=s(\beta^*_w+cf^*_w)$ and
$\tilde\beta_{uu}=s(\beta^*_{uu}+cf^*_{uu})$. So the conditions of
Definition~\ref{def:canard} hold for both systems or for neither. In
Lemma~\ref{lem:A} the two terms with $c$ cancel, so $\mathfrak B$ is
multiplied by $s$ and $A$ by $\sqrt s$.
\end{proof}

The map $\Phi$ preserves the fast fibres $y=\mathrm{const}$ and the
orientation of the fast variable, since $\partial u/\partial x=1/y>0$. At
corresponding points, $f=y^df^*$, $f_x=y^{d-1}f^*_u$ and
$\beta=-y^d(\beta^*+f^*)$. So $p$ satisfies $f=f_x=\beta=0$ if and only if
$\Phi(p)$ satisfies $f^*=f^*_u=\beta^*=0$, and then the chain rule gives
\[
  f_{xx}=y^{d-2}f^*_{uu},\qquad f_y=-y^{d-2}f^*_w,\qquad
  \beta_x=-y^{d-1}\beta^*_u,
\]
and, carried to third order, $\mathfrak B=y^{d-2}\mathfrak B^*$. Hence
$A(G^*)=\sqrt y\,A(G)$ at corresponding canard points. This is a second
proof of Theorem~\ref{thm:dual}, in the state space instead of in flux
coordinates. The fluxes of $G^*$ at $\Phi(x,y)$ are $y^{-d}v$. This is a
positive multiple of the fluxes of $G$, and the conditions~\eqref{eq:cone}
are invariant under positive scaling, which explains why the canard
conditions hold at ``the same fluxes''. Since $w=1/y$ reverses the slow
axis, $f_y$ and $\beta_x$ both change sign, and the two mechanisms are
exchanged.

The change of chart also explains the name. Pad every complex
$a\ce X+b\ce Y$ with $d-a-b$ vacancies V. The mass-action monomials become
homogeneous of degree $d$ in $(x,y,v)$. Those of $G$ are their restrictions
to $v=1$, and those of $G^*$ their restrictions to $y=1$, after V has been
renamed Y. So the duality is not a symmetry of one network: it relates two
different mass-action networks. At $\eps=0$ the system~\eqref{eq:chart} and
that of $G^*$ have the same layer problem, and their reduced flows on the
critical manifold differ by the positive factor $w$. The two networks
therefore have equivalent singular limits. For $\eps>0$ they are not
equivalent. Near the canard point, \eqref{eq:chart} is the system of $G^*$
with its slow equation multiplied by $w$, up to terms that vanish at the
canard point. These terms do not change the sign of $A$. They do change, by
amounts of the same order, the offsets of order $\eps$ of the Hopf point and
of the explosion. For instance, the distance $S\eps$ between them
(Remark~\ref{rem:S}) satisfies $S(G^*)=y\,S(G)$.

Chemically, the vacancy is a constant background pool. Every padded complex
has $d$ molecules, so every padded reaction conserves the total amount of X,
Y and V, just as reactions on a lattice or a catalytic surface conserve the
number of sites, with V as the empty site. In $G$ the concentration of V is
held constant and absorbed into the rate constants, like that of a solvent
or a buffered reservoir. In $G^*$ it is Y that is held constant, and V that
varies on the slow time scale. The duality therefore exchanges the slow
species with the background pool, and the chart $(x/y,1/y)$ measures X and
the pool in units of Y. In the planar autocatalator, the substrate Y is fed
from the pool by $0\to\ce Y$, and X decays into the pool by $\ce X\to0$. In
its dual, the decay product of X becomes the slow species
($\ce X+2\ce Y\to3\ce Y$), and the former substrate becomes part of the
constant background. So a substrate that X depletes becomes a product of X
that accumulates and removes X, and substrate depletion becomes
activator--inhibition.

\section{Elementary modes and the mode criterion}\label{sec:modes}

Suppose $K\neq\emptyset$. Then the closure $\bar K=\{v\ge0:Bv=0\}$ is a
pointed polyhedral cone, and $K$ is its relative interior. Its extreme rays
are the \emph{elementary modes}: nonnegative kernel vectors of minimal
support. A minimal support $T$ has $\rank B_T=|T|-1$, so every elementary
mode involves at most four reactions.

\begin{definition}
For an elementary mode $e$, let $p(e)=\big(n_{01}(e),m_{10}(e)\big)\in\R^2$
and $\Omega=\{(\varphi,\psi):\varphi\psi<0\}$.
\begin{itemize}[nosep]
  \item Two modes $e,e'$ (possibly equal) are \emph{sign-meeting} if the cone
        spanned by $p(e)$ and $p(e')$ meets $\Omega$.
  \item A mode $e$ \emph{carries a fold} if $n_{20}(e)\neq0$.
\end{itemize}
\end{definition}

\begin{theorem}[Mode criterion]\label{thm:modes}
A network $G$ with $K\neq\emptyset$ has a canard point iff its elementary modes
include a sign-meeting pair and a mode that carries a fold.
\end{theorem}

\begin{proof}
($\Rightarrow$) Every $v\in K$ is a combination $v=\sum_ic_ie_i$ of all the
elementary modes with all $c_i>0$, because $K$ is the relative interior of
$\bar K$. If $p(v)=\sum_ic_ip(e_i)\in\Omega$, then by Carathéodory's theorem
for cones in $\R^2$ the point $p(v)$ lies in the cone of at most two of the
$p(e_i)$. So two modes are sign-meeting. If $n_{20}(v)\neq0$, some mode has
$n_{20}(e_i)\neq0$.

($\Leftarrow$) Let $e,e'$ be sign-meeting, with $p(w_0)\in\Omega$ for some
$w_0=se+s'e'$, $s,s'\ge0$. Let $z$ carry a fold. For all
small $\gamma>0$, the point $w=w_0+\gamma z$ has $p(w)\in\Omega$, because
$\Omega$ is open, and $n_{20}(w)=n_{20}(w_0)+\gamma n_{20}(z)\neq0$ for all but
at most one value of $\gamma$. Finally, $v=w+\delta u$ with $u\in K$ and small
$\delta>0$ is strictly positive and still satisfies the open conditions
of~\eqref{eq:cone}.
\end{proof}

\begin{corollary}[At most three modes]\label{cor:three}
Every network with a canard point contains a subnetwork with a canard point
that is the union of the supports of at most three elementary modes: a
sign-meeting pair $e,e'$ and a mode $z$ carrying a fold. It therefore has at
most twelve reactions. If $e$ or $e'$ carries a fold, or $z$ is
sign-meeting with $e$ or with $e'$, then two of these modes suffice.
\end{corollary}

\begin{proof}
Let $W=\supp e\cup\supp e'\cup\supp z$. Then $e+e'+z$ is a strictly positive
flux of $W$, so Theorem~\ref{thm:modes} applies to $W$. If $e$ or $e'$
carries a fold, the pair suffices. If $z$ is sign-meeting with $e$ or with
$e'$, that pair suffices.
\end{proof}

\begin{corollary}[Minimal networks]\label{cor:minimal}
Call a canard network \emph{minimal} if no proper subnetwork has a canard
point. Then a network $U$ is minimal iff two conditions hold:
\begin{itemize}[nosep]
  \item it is the union of the supports of a sign-meeting pair and a mode
        carrying a fold, where the fold mode may belong to the pair;
  \item no other such choice among the elementary modes of $U$ has a strictly
        smaller union.
\end{itemize}
This can be decided exactly from the elementary modes supported in $U$.
\end{corollary}

\begin{figure}[ht]
\centering
\includegraphics[width=\textwidth]{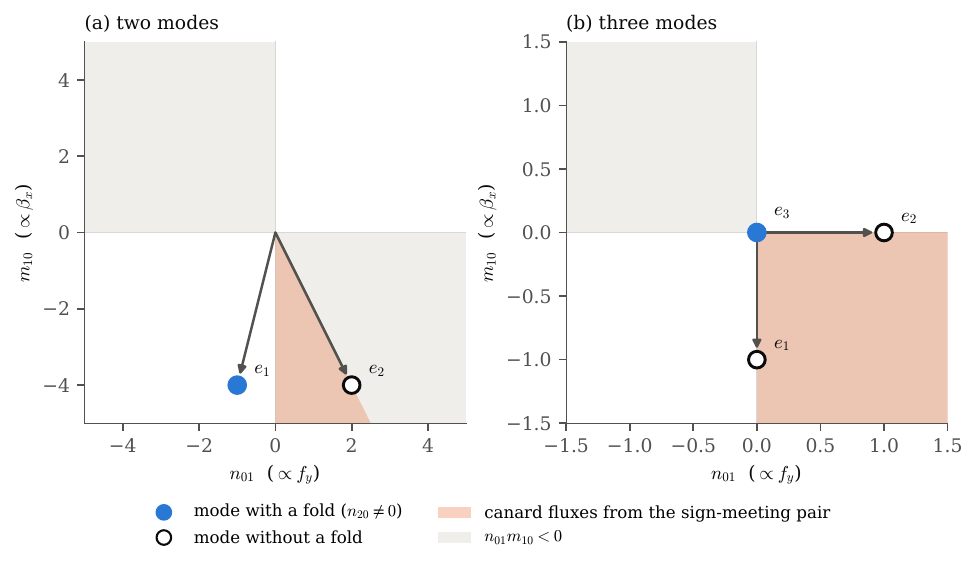}
\caption{The mode criterion in the $(n_{01},m_{10})$-plane. Each elementary
mode $e$ is drawn at $p(e)$; the shaded quadrants are $\Omega$. (a) The
five-reaction network of Example~\ref{ex:five}. Mode $e_2$ lies in $\Omega$ but
carries no fold, while $e_1$ carries the fold but lies outside $\Omega$.
Mixtures with $n_{01}>0$ are canard fluxes. (b) The six-reaction network of
Example~\ref{ex:three}. The sign-meeting pair $e_1,e_2$ carries no fold, and the
fold mode $e_3$ maps to the origin, so all three modes are needed.}
\label{fig:modes}
\end{figure}

All three cases occur (Table~\ref{tab:class}). A single mode can be
sign-meeting with itself and carry a fold. In the two-mode case, each mode
may fail on its own (Figure~\ref{fig:modes}(a)).
In the three-mode case, neither the pair nor any pair with the fold mode has a
canard point, but all three together do (Figure~\ref{fig:modes}(b)).

\begin{proposition}[Inheritance]\label{prop:inherit}
Let $W$ be a subnetwork of $G$, and suppose $W$ has a canard point at fluxes
$w$ with $\rank B_W=3$. Then $G$ has a canard point. If $Q(w)\neq0$, then $G$
has a canard point of the same criticality.
\end{proposition}

\begin{proof}
$B_W$ maps onto $\R^3$. Given small $u>0$ on the reactions of $G$ not in $W$,
solve $B_Wz=-B_{G\setminus W}u$ with $z=O(|u|)$. Then $v=(w+z,u)$ lies in the
canard cone of $G$, and the open conditions of~\eqref{eq:cone} and the sign
of $Q$ persist.
\end{proof}

\begin{remark}[New species]\label{rem:species}
Canard points, and the criticality of the singular Hopf bifurcation, also
persist when weakly coupled species are added. Let $W$
have a canard point with $\rank B_W=3$ and $Q\neq0$. Add species
$\ce Z_1,\dots,\ce Z_p$ with inflow and outflow $0\rightleftarrows\ce Z_i$
on the fast time scale, and reactions that couple them to X and Y with rate
constants $\delta k'$, in which every change of Y is slow. At $\delta=0$ the
concentrations $z_i$ relax exponentially to constants, independently of X and
Y. Let $\lambda$ be a rate constant of $W$ that unfolds the canard point
(Theorem~\ref{thm:planar}(d)). Extend the system by
$\dot\eps=\dot\delta=\dot\lambda=0$. Near the canard point it has a local
centre manifold of class $C^k$ for every finite $k$. For fixed parameters
this is a two-dimensional attracting locally invariant manifold. At
$\delta=0$ the equations for the $z_i$ are linear and decoupled, so
$z=z^*+O(\delta)$ on the manifold. On it the system reduces to
$\dot x=f+\delta f_1$, $\dot y=\eps(\beta+\delta\beta_1)$, with $f_1$ and
$\beta_1$ depending on $(x,y,\eps,\delta,\lambda)$; by the remark after
Lemma~\ref{lem:A}, the dependence of $f_1$ on $\eps$ does not matter. The conditions
of Definition~\ref{def:canard} are equations with a nonsingular Jacobian,
together with open inequalities, and $\lambda$ unfolds the canard point.
Hence, for small $\delta$, the reduced system has a canard point at a nearby
value of $\lambda$, and its coefficient $A$ is close to that of $W$. In
particular, the criticality is the same, and since the manifold attracts,
so is that of the singular Hopf bifurcation of the full system. The
relaxation oscillations and the large canard cycles need an invariant
manifold along the whole critical manifold, which this local argument does
not provide.
\end{remark}

\section{Classification for molecularity at most three}\label{sec:class}

The complexes $a\ce X+b\ce Y$ with $a+b\le3$ give $90$ reactions. By
Corollary~\ref{cor:minimal}, the minimal canard networks can be listed from
the elementary modes of the $3\times90$ balance matrix, without enumerating
subsets of reactions. The enumeration proceeds in four steps, followed by the
labelling.

\paragraph{Modes.}
There are $15\,804$ elementary modes: $130$ of support $2$, $1\,368$ of
support $3$ and $14\,306$ of support $4$. Of these, $4\,982$ (of support $4$)
and $16$ (of support $3$ and rank $2$, hence of codimension zero) are
sign-meeting with themselves and carry a fold. They are the minimal canard
networks carried by a single mode.

\paragraph{Pairs of modes.}
Among the remaining modes we consider the pairs in which at least one member
carries a fold. There are $55\,437\,890$ such pairs, and the sign-meeting
ones give $18\,370\,633$ distinct unions of supports.

\paragraph{Triples of modes.}
A network that needs three modes contains a sign-meeting pair in which
neither member carries a fold and whose union contains no fold mode; there
are $419\,332$ such pairs. To each we add a fold mode that is sign-meeting
with neither member and whose new reactions form an inclusion-minimal set
among all such fold modes. Both conditions are necessary for minimality, and
together they give $3\,917\,188$ distinct unions.

\paragraph{Minimality.}
By Corollary~\ref{cor:minimal}, a union $U$ is minimal if and only if no
canard set of one, two or three modes inside $U$ has a strictly smaller
union. This test is exact, and it selects the minimal networks among the
candidates of the previous steps.

\paragraph{Criticality.}
Each minimal network receives one of four labels. It is \emph{super only}
if $A\le0$ on all admissible fluxes, with $A<0$ on an open dense subset, and
\emph{sub only} in the analogous case; it is labelled \emph{both} if both
signs occur, and $A\equiv0$ if $Q$ vanishes on the cone. For five reactions
with $\rank B=3$, the cone is $\{e_1+te_2:t>0\}$, and $Q$ is a cubic in $t$
with integer coefficients; its signs on the admissible set are found exactly,
by rational evaluation between all breakpoints and roots. Otherwise, the
admissible set is a union of at most four open convex polyhedral cones, one
for each quadrant of $(n_{01},m_{10})$ and each sign of $n_{20}$. On each
piece we write $Q$ as a cubic form in the coefficients $\mu\ge0$ of its
extreme rays. If all coefficients of these forms have the same sign on every
piece, then $Q$ has that sign wherever it is nonzero; this is Pólya's
certificate of degree $0$. Zeros on the walls between pieces are allowed,
since they are isolated or lower-dimensional. The label \emph{both} is
certified by exact rational fluxes of each sign, and $A\equiv0$ by the
identity $Q\equiv0$. Every network received an exact label.

\begin{table}[ht]
\centering
\small
\begin{tabular}{@{}rrrrrrr@{}}
\toprule
reactions & modes & super only & sub only & both & $A\equiv0$ & total\\
\midrule
\multicolumn{7}{@{}l}{\emph{codimension one} ($\rank B=3$)}\\
$4$ & $1$ & $2\,460$ & $2\,314$ & -- & $208$ & $4\,982$\\
$5$ & $2$ & $19\,142$ & $34\,418$ & $19\,266$ & $1\,132$ & $73\,958$\\
$6$ & $2$ & $188$ & $1\,070$ & $20\,974$ & $30$ & $22\,262$\\
$6$ & $3$ & $2\,708$ & $3\,714$ & $31\,080$ & $18$ & $37\,520$\\
\midrule
\multicolumn{7}{@{}l}{\emph{codimension zero} ($\rank B=2$)}\\
$3$ & $1$ & $4$ & $12$ & -- & -- & $16$\\
$4$ & $2$ & $16$ & -- & -- & -- & $16$\\
\midrule
total & & $24\,518$ & $41\,528$ & $71\,320$ & $1\,388$ & $138\,754$\\
\bottomrule
\end{tabular}
\caption{All minimal planar canard networks with complexes of molecularity at
most three, X fast and Y slow, by number of reactions and number of modes.
For six reactions, two rows are needed: networks carried by a pair of modes,
and networks that need three. All labels are exact. The list is the file
\texttt{data/minimal\_canard\_networks.tsv} of the data record
\citep{CanardPlanarData2026}.}
\label{tab:class}
\end{table}

\paragraph{Checks.}
Five independent checks support the classification, and a sixth computation
tests the labels for finite $\eps$.

\paragraph{Two routes.}
An independent enumeration of all $4\,773\,516$ five-reaction networks that
contain $2\ce X\to3\ce X$ or $2\ce X+\ce Y\to3\ce X$ yields the same
$73\,958$ minimal networks.

\paragraph{Minimality by linear programming.}
This check tests minimality without elementary modes, without
Carathéodory's theorem and without Theorem~\ref{thm:modes}. It uses only the
balance matrix and the three linear forms in~\eqref{eq:cone}. Condition
\eqref{eq:cone} is not convex, because of the inequation $n_{20}\neq0$ and
the product $n_{01}m_{10}<0$. It is, however, the union of four convex
conditions, one for each choice of signs $s=(s_1,s_2,s_3)\in\{\pm1\}^3$ with
$s_1s_2=-1$:
\[
  v>0,\qquad Bv=0,\qquad s_1n_{01}(v)>0,\qquad s_2m_{10}(v)>0,\qquad
  s_3n_{20}(v)>0 .
\]
Each defines an open polyhedral cone. Since all these conditions are
invariant under $v\mapsto\lambda v$ with $\lambda>0$, a solution can be
rescaled so that every strict inequality holds with margin at least $1$;
conversely, a solution with these margins satisfies the strict
inequalities. Hence a network $W$ has a canard point for some rate constants
if and only if at least one of the four polyhedra
\[
  P_s(W)=\bigl\{v\in\R^{W}:\ v\ge\mathbf 1,\ Bv=0,\ s_1n_{01}(v)\ge1,\
  s_2m_{10}(v)\ge1,\ s_3n_{20}(v)\ge1\bigr\}
\]
is nonempty. Each is a linear feasibility problem with $|W|+3$ inequalities
and three equations. We solve it in exact rational arithmetic (polymake,
through Oscar), so that no numerical tolerance enters.

The network $W$ is minimal if and only if some $P_s(W)$ is nonempty and
every $P_s(W')$ is empty for every proper subnetwork $W'$. Subnetworks with
at most two reactions never have a canard point, so they need not be tested.
Indeed, if $v>0$ satisfies $n_{00}=n_{10}=0$ on two reactions, then either
$N_1=N_2=0$, or $N_1v_1=-N_2v_2\neq0$ and $a_1=a_2$. In both cases
$n_{20}=a_1(a_1-1)N_1v_1+a_2(a_2-1)N_2v_2=0$, and the case of one reaction
is the same. All other proper subnetworks are tested, not only those with one
reaction fewer. By Proposition~\ref{prop:inherit}, a canard point of a
subnetwork with $\rank B=3$ is inherited by every larger network. A canard
point of rank $2$, however, persists only in larger networks whose cone is
nonempty (Theorem~\ref{thm:modes}). A network whose subnetworks with one
reaction fewer all fail may therefore still contain a canard network of
codimension zero. For a network of six reactions the test involves the $41$
subnetworks with three to five reactions, each with four feasibility
problems.

We applied the test to random samples (with a fixed seed) of $200$ networks
from each of the four classes of codimension one in Table~\ref{tab:class},
and to all $32$ networks of codimension zero: $832$ networks in total. All
$832$ were confirmed to be minimal, in about nine minutes. The test checks
the minimality of the sampled networks; that the list is complete rests on
Theorem~\ref{thm:modes}, and for five reactions it is confirmed by the
independent enumeration above.

\paragraph{Duality.}
Theorem~\ref{thm:dual} maps the list onto itself with the same labels, for
all $138\,754$ networks, and no network is self-dual. Since each exact label
was certified separately for each network, this is a test of
Theorem~\ref{thm:dual}. For $d=3$ the duality exchanges
$2\ce X+\ce Y\to3\ce X$ with $2\ce X\to3\ce X$. Accordingly, $69\,220$
minimal networks contain only the first, $69\,220$ only the second, with the
same distribution of labels, and $314$ contain both.

\paragraph{Direct computation.}
The labels of the four-reaction networks agree with an independent
enumeration based on Lemma~\ref{lem:A}. For $25$ sampled four-reaction networks and several
larger ones, we confirmed the canard point and the sign of $A$ directly from
the rate constants, using the first Lyapunov coefficient of the full system
as $\eps\to0$.

\paragraph{Independent verification.}
A separate program re-derives every entry of the list from the reactions
alone. It is written in Python with exact integer and rational arithmetic,
and uses none of the Julia code, Oscar or polymake. For each of the
$138\,754$ networks it recomputes the elementary modes. By
Theorem~\ref{thm:modes} it then checks that the network has a canard point
and that none of its proper subnetworks has one, and it checks the number of
modes and the rank. Next it recomputes the polyhedral pieces of the
admissible set and their extreme rays, and certifies each label. The labels
\emph{super} and \emph{sub} are certified by Pólya coefficients of one sign
on every piece ($65\,780$ networks). Where these coefficients have mixed
signs on a two-dimensional piece, the program isolates the real roots of the
cubic exactly with Sturm sequences ($266$ networks). The label \emph{both} is
certified by two rational admissible fluxes at which $Q$ has opposite signs
($71\,320$), and $A\equiv0$ by the identity $Q\equiv0$ on the flux cone
($1\,388$). All $138\,754$ entries pass, the list is closed under the
duality, and the counts of Table~\ref{tab:class} are reproduced. The witness
fluxes are part of the data record, so that every label \emph{both} can be
checked by hand.

\paragraph{Criticality at finite $\eps$.}
We also tested how well the labels predict the criticality of the Hopf
bifurcation for finite $\eps$. We took $50$ random networks from each row of
codimension one of Table~\ref{tab:class} and dropped two with $A\equiv0$.
For each, we chose a random admissible flux $v$, placed the canard point at
$(1,1)$, and unfolded it by the rate constant with the largest effect on
$\beta$ at the fold. We then continued the Hopf point in $\eps$ from
$10^{-5}$ to $0.2$ and computed its first Lyapunov coefficient $l_1$. Wherever
the Hopf point exists, the sign of $l_1$ agrees with that of $Q(v)$ in $188$
of $189$ cases at $\eps=10^{-3}$, in $135$ of $146$ at $\eps=10^{-2}$ and in
$74$ of $86$ at $\eps=0.1$. Of the $24$ networks with a disagreement at some
$\eps$, $20$ are labelled \emph{both}; their fluxes lie near degenerate canard
points, where the corrections of order $\sqrt\eps$ to $A$ decide. As $\eps$
grows, the Hopf branch ends for $51$ of the $198$ networks before
$\eps=10^{-2}$, and for $104$ before $\eps=0.1$; in $96$ of these $104$ the
equilibrium becomes a saddle, so the branch ends at a double zero
eigenvalue. The estimate $\eps_{\mathrm{BT}}$ of Remark~\ref{rem:BT} locates
this end when it is small. For the $68$ networks with
$\eps_{\mathrm{BT}}<0.03$, the observed end lies within a factor of two of
$\eps_{\mathrm{BT}}$ in $65$ cases, with median ratio $0.99$. For larger
$\eps_{\mathrm{BT}}$ the ratio ranges from $0.15$ to $2$, and $70$ branches
persist up to $\eps=0.2$. For $49$ of the $198$ sampled fluxes,
$\eps_{\mathrm{BT}}<10^{-2}$, so the time scales must be separated by about
two orders of magnitude or more, and $\eps_{\mathrm{BT}}$ says in advance by
how much.

\begin{example}[A minimal five-reaction network]\label{ex:five}
$\{0\to\ce Y,\ 0\to\ce X,\ \ce X+\ce Y\to0,\ 2\ce X\to0,\ 2\ce X+\ce Y\to3\ce X\}$.
Here $K=\{2v_4<v_5\}$, with $n_{01}=4v_4-v_5$, $m_{10}=4v_4-4v_5$ and
$n_{20}=2v_5-4v_4$. The canard condition is $2v_4<v_5<4v_4$, and it fails at
both extreme rays: at $v_4=0$ the four-reaction subnetwork has the saddle
sign, and at $v_5=2v_4$ the three-reaction subnetwork has $n_{20}=0$.
At $v=(5,1,2,1,3)$ and $(x,y)=(1,1)$, $Q=-56$, and the rate-constant
computation gives $A=-9.90$.
\end{example}

\begin{figure}[ht]
\centering
\includegraphics[width=\textwidth]{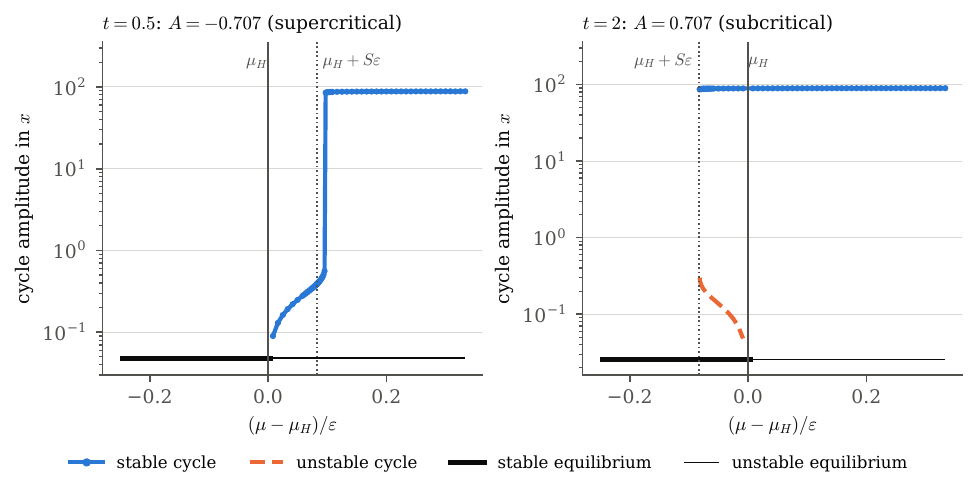}
\caption{Canard explosions of the network of Example~\ref{ex:both} for
$\eps=0.01$. The horizontal axis is $(\mu-\mu_H)/\eps$, where $\mu_H$ is the
exact singular Hopf point. The amplitude in $x$ of the stable and unstable
limit cycles is shown on a log scale. Stable cycles were found by forward
integration, unstable ones by backward integration (they attract in reversed
time), and the equilibrium is drawn along the bottom. Left, $t=1/2$
($A<0$): stable small cycles grow from $\mu_H$ and explode into relaxation
oscillations. Right, $t=2$ ($A>0$): unstable small cycles coexist with the
stable equilibrium and a stable relaxation cycle, so the system is bistable.
In both cases the explosion is close to the first-order prediction
$\mu_H+S\eps$ of Remark~\ref{rem:S} ($S=\pm1/12$, dotted); the difference is
of the expected order $\eps^{3/2}$ or smaller. The canard cycles themselves lie in an exponentially
small $\mu$-interval and appear as the vertical part of the branch.}
\label{fig:explosion}
\end{figure}

\begin{example}[A degenerate canard point]\label{ex:both}
$\{0\to\ce X,\ \ce X\to\ce Y,\ \ce X+\ce Y\to\ce Y,\ \ce X+\ce Y\to2\ce X,\ 2\ce X+\ce Y\to3\ce X\}$.
Its two modes are $e_1$ on $\{\ce X\to\ce Y,\ \ce X+\ce Y\to2\ce X\}$ and
$e_2$ on
$\{0\to\ce X,\ \ce X\to\ce Y,\ \ce X+\ce Y\to\ce Y,\ 2\ce X+\ce Y\to3\ce X\}$,
each with all fluxes equal to $1$. Each fails on its own: $e_1$ has
$n_{20}=m_{10}=0$, and $e_2$ has $f_y=0$.
Along $v=e_1+te_2$ we have $n_{01}=1$, $m_{10}=-t$ and $n_{20}=2t$, and
\[
  Q=2t^2(t-1).
\]
The explosion is supercritical for $t<1$ and subcritical for $t>1$. At
$t=1$, i.e.\ at rate constants $(1,2,1,1,1)$ and $(x,y)=(1,1)$, the canard
point is degenerate with $A=0$. This is a codimension-two point in
rate-constant space.

Figure~\ref{fig:explosion} shows the explosion on either side of $t=1$, using
as unfolding parameter $\mu$ the rate constant of $0\to\ce X$. The other rate
constants are $(1+t,\ t,\ 1,\ t)$, and the canard point is at $\mu=t$.
\end{example}

\begin{example}[A network that needs three modes]\label{ex:three}
$\{0\to\ce Y,\ 0\to\ce X,\ \ce X\to0,\ \ce X+\ce Y\to\ce X,\ \ce X+\ce Y\to2\ce X+\ce Y,\ 2\ce X\to3\ce X\}$.
It has three elementary modes: $e_1$ on
$\{0\to\ce Y,\ \ce X+\ce Y\to\ce X\}$, with $p(e_1)=(0,-1)$ and no fold;
$e_2$ on $\{\ce X\to0,\ \ce X+\ce Y\to2\ce X+\ce Y\}$, with $p(e_2)=(1,0)$
and no fold; and $e_3$ on $\{0\to\ce X,\ \ce X\to0,\ 2\ce X\to3\ce X\}$, with
fluxes $1,2,1$, $p(e_3)=(0,0)$ and a fold.
Only $e_1,e_2$ are sign-meeting, and neither carries a fold. So no
subnetwork carried by one or two modes has a canard point, but the whole
network does. For $v=ae_1+be_2+ce_3$,
\[
  n_{01}=b,\qquad m_{10}=-a,\qquad n_{20}=2c,\qquad Q=2ac\,(2c-b),
\]
so every positive mixture is admissible (substrate depletion). The explosion
is supercritical for $b>2c$ and subcritical for $b<2c$. For instance, at
$a=b=c=1$ the rate constants are $(1,1,3,1,1,1)$ at $(x,y)=(1,1)$: then
$f=1-3x+xy+x^2$ and $\beta=1-xy$ satisfy $f=f_x=\beta=0$, and $Q=2>0$.
\end{example}

\section{Realisation by bimolecular steps}\label{sec:bimol}

By Corollary~\ref{cor:coop}, a planar canard network needs a step with three
molecules. In chemistry such steps are effective descriptions, usually of
bimolecular steps through fast intermediates. We show that canard points
survive this resolution, and that the criticality changes by an explicit
term.

Suppose that X and Y are also stored in further species by fast reversible
reactions, so that at the equilibrium of these reactions the total amounts of
X and Y are $\varphi(x)$ and $\psi(y)$, with $\varphi'>0$ and $\psi'>0$.
More precisely, we assume the following:
\begin{itemize}[nosep]
  \item the storage reactions conserve the total amounts of X and of Y, and
        store X and Y separately;
  \item at fixed totals their equilibrium is unique and normally attracting;
  \item the other reactions keep their mass-action rates in the free
        concentrations $x$ and $y$, and change the totals by $N_r$ and
        $M_r$;
  \item every change of the total of Y by these reactions is slow.
\end{itemize}
On the critical manifold of the storage reactions, the totals then obey
$\tfrac{d}{dt}\varphi(x)=f$ and $\tfrac{d}{dt}\psi(y)=\eps\beta$, that is,
\begin{equation}\label{eq:storage}
  \dot x=\frac{f(x,y)}{\varphi'(x)},\qquad \dot y=\eps\,\frac{\beta(x,y)}{\psi'(y)} .
\end{equation}
This is the rapid-buffer approximation of \citet{WagnerKeizer1994}; what is
new below is its effect on the criticality.

\begin{proposition}[Fast storage]\label{prop:storage}
The system~\eqref{eq:storage} has the same canard points as~\eqref{eq:sys},
at the same fluxes. At a canard point, $\sign A=\sign Q_\theta(v)$, where
\[
  Q_\theta=Q+2\theta\,n_{01}m_{10}n_{20},\qquad
  \theta=\frac{x\varphi''(x)}{\varphi'(x)} .
\]
In particular, the storage of Y does not affect the criticality.
\end{proposition}

\begin{proof}
Let $\tilde f=rf$ and $\tilde\beta=s\beta$, with $r=1/\varphi'(x)$ and
$s=1/\psi'(y)$, and let $f=f_x=\beta=0$ at $p$. At $p$,
\[
  \tilde f_x=0,\quad \tilde f_{xx}=rf_{xx},\quad \tilde f_y=rf_y,\quad
  \tilde f_{xy}=rf_{xy}+r'f_y,\quad \tilde f_{xxx}=rf_{xxx}+3r'f_{xx},
\]
and $\tilde\beta_x=s\beta_x$, $\tilde\beta_y=s\beta_y$,
$\tilde\beta_{xx}=s\beta_{xx}$. So the conditions of
Definition~\ref{def:canard} hold for both systems or for neither, and
Lemma~\ref{lem:A} gives
\[
  \tilde{\mathfrak B}=s\Bigl(\mathfrak B-2\,\frac{r'}{r}\,\frac{f_y\beta_x}{f_{xx}}\Bigr).
\]
Multiplying by $x^4y\,f_{xx}^2>0$ and using~\eqref{eq:moments} gives
$sQ_\theta(v)$, with $\theta=-xr'/r=x\varphi''/\varphi'$.
\end{proof}

Since $n_{01}m_{10}<0$, the correction has the sign of $-\theta n_{20}$.
Storage that is cooperative in X ($\varphi''>0$ at the canard point,
$\theta>0$) pushes the explosion towards supercritical where $f_{xx}>0$, and
towards subcritical where $f_{xx}<0$; a saturable buffer ($\varphi''<0$,
$\theta<0$) does the opposite. For a fast dimerisation
$2\ce X\rightleftarrows\ce D$ with equilibrium constant $K$, we have
$\varphi(x)=x+2Kx^2$ and
\[
  \theta=\frac{4Kx}{1+4Kx}\in(0,1),
\]
where $x$ is the X-coordinate of the canard point. So $\theta$ increases with
the dimensionless strength $4Kx$ of the dimerisation. Dimers give a
bimolecular realisation of every network of codimension one in the
classification.

\begin{corollary}[Bimolecular realisation]\label{cor:bimol}
Let all complexes of $G$ have molecularity at most three. Replace the
complexes $3\ce X$, $2\ce X+\ce Y$, $\ce X+2\ce Y$ and $3\ce Y$ by
$\ce D+\ce X$, $\ce D+\ce Y$, $\ce X+\ce E$ and $\ce E+\ce Y$, where D and E
are new species. Add the reactions $2\ce X\rightleftarrows\ce D$ and
$2\ce Y\rightleftarrows\ce E$, with equilibrium constants $K$ and $K'$, on a
time scale $\delta$ faster than that of X, and divide each rate constant by
$K$ or $K'$ for each D or E in its reactant. In the resulting network
$G^{\mathrm{bi}}$, every change of Y and of E by a reaction other than
$2\ce Y\rightleftarrows\ce E$ is slow, as in~\eqref{eq:sys}, and all
complexes have at most two molecules. Suppose that $G$ has a
canard point at fluxes $v$, with $\rank B=3$ and $Q_\theta(v)\neq0$ for
$\theta=4Kx/(1+4Kx)$. Then, for all small $\delta>0$, $G^{\mathrm{bi}}$ has a
two-dimensional locally invariant manifold on which, for suitable rate
constants close to the given ones, it has a canard point of the criticality
given by $\sign Q_\theta(v)$.
\end{corollary}

\begin{proof}
At $\delta=0$, the equilibria $d=Kx^2$ and $e=K'y^2$ of the dimerisations
form a normally attracting manifold. On it the rates of $G^{\mathrm{bi}}$ are
those of $G$, and the total amounts $x+2d$ and $y+2e$ are $\varphi(x)$ and
$\psi(y)=y+2K'y^2$. Every change of $y+2e$ is slow, because the dimerisation
of Y conserves it and the other reactions change Y and E slowly. By the
theorem of \citet{Fenichel1979}, the manifold persists for small $\delta>0$
as an attracting locally invariant manifold, of class $C^k$ in
$(\delta,\eps)$ for every finite $k$. In the coordinates $(x+2d,\,y+2e)$ the
system on it is a perturbation of order $\delta$ of~\eqref{eq:storage}
written in $(\varphi(x),\psi(y))$, of the same slow--fast form. A rate
constant unfolds the canard point of~\eqref{eq:storage}
(Theorem~\ref{thm:planar}(d)); the factors $r$ and $s$ do not affect this.
As in Remark~\ref{rem:species}, the canard point therefore persists for
small $\delta$, with $A$ close to that of~\eqref{eq:storage}, whose sign is
$\sign Q_\theta(v)$ by Proposition~\ref{prop:storage}. Since the manifold
attracts, the singular Hopf bifurcation of $G^{\mathrm{bi}}$ has this
criticality.
\end{proof}

For a given $K$, moving the canard point along the X-axis at fixed fluxes
realises every $\theta\in(0,1)$. So $\theta$ is a property of the canard
point rather than of the network, and the labels below are labels at a
given strength $4Kx$. For the four-reaction networks of
Table~\ref{tab:four}, the flux is unique up to scaling and $Q_\theta$ is
linear in $\theta$. So a label holds for every strength of the dimerisation
if and only if $Q_1$ is zero or has the sign of $Q$. Of the $4\,774$ networks
with $A\not\equiv0$, $4\,416$ keep their label for every $\theta\in[0,1)$,
and $358$ change it at $\theta^*=-Q/(2n_{01}m_{10}n_{20})\in(0,1)$: $256$
from subcritical to supercritical and $102$ the other way. As networks,
the realisations $G^{\mathrm{bi}}$ of these $358$ therefore admit both
criticalities. The $208$
networks with $A\equiv0$ become nondegenerate: $156$ supercritical and $52$
subcritical, according to the sign of $f_{xx}$.

For the whole classification, we computed the labels of $Q_1$, the limit of
strong dimerisation, exactly and in the same way as those of $Q$ (see
\emph{Independent verification} in Section~\ref{sec:class}); none remained
undetermined. Since $Q_\theta=(1-\theta)Q+\theta Q_1$, a one-sided label
holds for every $\theta\in[0,1)$ if and only if $Q_1$ does not take the
opposite sign. Of the $66\,014$ networks of codimension one with a one-sided
label, $57\,796$ keep it for every strength of the dimerisation. For the
other $8\,218$, the opposite sign appears once $\theta$ exceeds a threshold,
which can be $0$. For instance,
$\{0\to\ce Y,\ \ce Y\to\ce X,\ \ce X\to0,\ 2\ce X\to3\ce X,\ 2\ce X+\ce Y\to2\ce X\}$
is labelled \emph{sub}. Along its admissible fluxes $v=(1,s,2s,s,1-s)$,
$0<s<1$, the sign of $Q_\theta$ changes at $\theta^*=s/(2(1-s))$, which
tends to $0$ with $s$, so arbitrarily weak dimerisation creates supercritical
canard points. In the limit $\theta\to1$, $7\,450$ of the $8\,218$ become
\emph{both} and $768$ change their one-sided label. Of the $71\,320$
networks labelled \emph{both}, $62\,406$ still admit both criticalities in
the limit, and $8\,914$ become one-sided. The $1\,388$ networks with
$A\equiv0$ become nondegenerate, $988$ supercritical and $400$ subcritical.

We checked Corollary~\ref{cor:bimol} on three four-reaction networks. For
each, we computed the first Lyapunov coefficient $l_1$ of the full network
$G^{\mathrm{bi}}$ on four species, at its Hopf point located by real root
isolation on the Hurwitz condition, for $(\eps,\delta)=(10^{-4},10^{-3})$
and $(10^{-6},10^{-5})$, so in the joint limit with $\delta=10\eps$. For
$\{0\to\ce X,\ \ce X\to\ce Y,\ \ce X+\ce Y\to0,\ 2\ce X+\ce Y\to3\ce X\}$, with
$Q=8$ and $\theta^*=1/2$, $l_1$ changes from positive to negative as $K$
increases. For
$\{0\to\ce Y,\ \ce X\to0,\ 2\ce X+\ce Y\to3\ce X,\ 3\ce X\to2\ce X+\ce Y\}$,
with $Q=-4$ and $\theta^*=1/2$, it changes from negative to positive. In
both cases the change occurs within $2.4\cdot10^{-4}$ of $\theta^*$ for the
first pair $(\eps,\delta)$ and within $4\cdot10^{-5}$ for the second. The
network
$\{\ce Y\to\ce X,\ \ce X\to0,\ \ce X+2\ce Y\to3\ce Y,\ 2\ce X\to3\ce X\}$
needs both dimers and has $\theta^*=-2$; it stays supercritical for the five
tested pairs $(K,K')$, with $K$ from $1/100$ to $2$ and $K'$ from $1/100$ to
$100$, although $K'$ changes the size of $l_1$. In all $30$ evaluations away from $\theta^*$, the sign of $l_1$
is that of $Q_\theta$. Between the two pairs $(\eps,\delta)$, $\sqrt\eps\,l_1$
changes by at most six per cent, the largest change occurring for the
strongest dimerisation of Y.

\section{Discussion}\label{sec:disc}

\paragraph{Minimal networks.}
For two species, the question of which mass-action networks admit a canard
explosion, and of which criticality, is answered by the minimal networks.
These have at most twelve reactions, whatever the molecularity. The answer is
governed by four objects:
\begin{itemize}[nosep]
  \item the canard cone;
  \item the two linear forms $n_{01}$ and $m_{10}$, which carry the
        mechanism;
  \item the fold form $n_{20}$;
  \item the cubic $Q$.
\end{itemize}
In general, minimal canard networks are \emph{not} the four-reaction ones.
Mixtures of two or three elementary modes produce further minimal networks.
The general bound is twelve reactions, but in the class of molecularity at
most three none has more than six. By Theorem~\ref{thm:modes}, every network
$G\supseteq W$ with $K_G\neq\emptyset$ that contains a minimal canard
network $W$ has canard points. If $\rank B_W=3$,
Proposition~\ref{prop:inherit} adds that $K_G\neq\emptyset$ holds
automatically and that the criticality is kept. In the sense of inheritance
results such as \citet{BanajiBorosHofbauer2025}, the minimal networks are the
candidates for ``atoms'' of canard behaviour.

\paragraph{Criticality.}
The sign of the cubic $Q$ is constant on the four-reaction cones, but it varies on
larger ones. Most minimal networks of six reactions admit both
criticalities ($52\,054$ of $59\,782$), separated by degenerate canard points
as in Examples~\ref{ex:both} and~\ref{ex:three}. These networks are natural places to look for
codimension-two canard phenomena.

\paragraph{Design rules.}
In chemical terms, the results give four rules for a planar network with a
canard explosion, X fast and Y slow. First, X needs cooperative
autocatalysis: a step that produces X from a reactant with at least two X
(Corollary~\ref{cor:coop}), with orders of production and consumption that
alternate (Proposition~\ref{prop:fold}). This is what allows the balance of
X, for frozen Y, to have a fold. Second, X must also be supplied by a step
without X in its reactant, or removed by a step with three X in its reactant
(Remark~\ref{rem:orient}). Third, the slow feedback must have one of two
signs: either Y inhibits the net production of X while X promotes that of Y
(activator--inhibitor), or Y feeds X while X consumes Y (substrate
depletion). Fourth, these ingredients need at most three pathways: one that
balances X across the fold and a pair that supplies the feedback
(Theorem~\ref{thm:modes}). A single rate constant then places the canard
point (Theorem~\ref{thm:planar}(d)).

The criticality is not fixed by the reactions alone, since most networks of
six reactions admit both. Where it is fixed, as for four reactions, the
orientation of the fold is a good guide. Among the four-reaction networks of
Table~\ref{tab:four} with $A\not\equiv0$, $3\,010$ have $f_{xx}>0$ at the
fold, so that the supply of X outweighs its third-order removal, and
$2\,264$ of these are supercritical. Of the $1\,764$ with $f_{xx}<0$,
$1\,568$ are subcritical. Table~\ref{tab:repr}
lists representative networks: the two classical planar models in the list,
their duals, two networks whose reactants have at most two molecules, and the
examples of Section~\ref{sec:class}.

\begin{table}[ht]
\centering
\footnotesize
\begin{tabular}{@{}rlcc@{}}
\toprule
& network & mechanism & criticality\\
\midrule
1 & $0\to\ce Y,\ \ce Y\to\ce X,\ \ce X\to0,\ 2\ce X+\ce Y\to3\ce X$
  & SD & super\\
2 & $0\to\ce X,\ \ce X\to0,\ \ce X\to\ce Y,\ 2\ce X+\ce Y\to3\ce X$
  & SD & super\\
3 & $3\ce Y\to2\ce Y,\ 2\ce Y\to\ce X+2\ce Y,\ \ce X+2\ce Y\to3\ce Y,\ 2\ce X\to3\ce X$
  & AI & super\\
4 & $3\ce Y\to\ce X+2\ce Y,\ \ce X+2\ce Y\to3\ce Y,\ \ce X+2\ce Y\to2\ce Y,\ 2\ce X\to3\ce X$
  & AI & super\\
5 & $0\to\ce X,\ 2\ce X\to\ce Y,\ \ce X+\ce Y\to0,\ 2\ce X\to3\ce X$
  & AI & super\\
6 & $0\to\ce X+\ce Y,\ \ce X\to0,\ \ce X+\ce Y\to2\ce X,\ 2\ce X\to3\ce X$
  & SD & sub\\
7 & $0\to\ce Y,\ 0\to\ce X,\ \ce X+\ce Y\to0,\ 2\ce X\to0,\ 2\ce X+\ce Y\to3\ce X$
  & SD & super\\
8 & $0\to\ce X,\ \ce X\to\ce Y,\ \ce X+\ce Y\to\ce Y,\ \ce X+\ce Y\to2\ce X,\ 2\ce X+\ce Y\to3\ce X$
  & SD & both\\
9 & $0\to\ce Y,\ 0\to\ce X,\ \ce X\to0,\ \ce X+\ce Y\to\ce X,\ \ce X+\ce Y\to2\ce X+\ce Y,\ 2\ce X\to3\ce X$
  & SD & both\\
10 & $0\to\ce X+\ce Y,\ \ce X\to0,\ 2\ce X+\ce Y\to3\ce X$
  & SD & super\\
\bottomrule
\end{tabular}
\caption{Representative minimal canard networks (X fast, Y slow). Rows 1
and 2: the planar autocatalator (Example~\ref{ex:auto}) and the
Brusselator. Rows 3 and 4: their duals (Theorem~\ref{thm:dual}). Rows 5
and 6: all reactants have at most two molecules. Rows 7--9: the networks of
Examples~\ref{ex:five}, \ref{ex:both} and~\ref{ex:three}. Row 10: a
codimension-zero network, without an explosion. Mechanism:
activator--inhibitor (AI) or substrate depletion (SD) at the admissible
fluxes. Criticality: the exact label of Table~\ref{tab:class}.}
\label{tab:repr}
\end{table}

\paragraph{Robust networks.}
For experiments, the most useful networks are those whose behaviour does not
depend on fine details. Table~\ref{tab:robust} lists five networks chosen by
three criteria. All of them have four reactions. Their flux is then unique up
to scaling, so by Corollary~\ref{cor:four} the canard point exists and has
the same criticality for all rate constants.
\begin{itemize}[nosep]
  \item \emph{Storage.} $Q$ and $Q_1=Q+2n_{01}m_{10}n_{20}$ have the same
        sign. By Proposition~\ref{prop:storage}, the criticality is then the
        same under fast storage of X of any strength $\theta\in[0,1)$, and in
        particular under every fast dimerisation of X.
  \item \emph{Time scales.} Let $\eta=-n_{01}m_{10}/m_{01}^2$. This is the
        estimate $\eps_{\mathrm{BT}}$ of Remark~\ref{rem:BT} for a canard
        point with $x=y$. In general $\eps_{\mathrm{BT}}=(y/x)\,\eta$, where
        $\eps x/y$ is the ratio of the relative rates $\dot y/y$ and
        $\dot x/x$ for fluxes of comparable size. For $m_{01}=0$ the estimate
        predicts no end of the Hopf branch. As in Section~\ref{sec:class}, we
        placed the canard point at $(1,1)$ and continued the Hopf point in
        $\eps$, here up to $\eps=1$.
  \item \emph{Molecularity.} In rows 1, 2, 4 and 5 no reactant has more than
        two molecules, so all rates are at most quadratic and no step has to
        be resolved. Row 3, the autocatalator, has a trimolecular step. By
        the first criterion, its bimolecular realisation through a fast
        dimer (Corollary~\ref{cor:bimol}) is supercritical for every $K$.
\end{itemize}
Wherever the Hopf point exists, the sign of its first Lyapunov coefficient
agrees with the label at all values of $\eps$ tested. In rows 1 and 2 the
Hopf bifurcation stays supercritical up to $\eps=1$, that is, without any
separation of time scales; a sharp explosion, of course, still needs small
$\eps$. In row 5 the Hopf branch ends at $\eps=0.104$, close to
$\eta=1/9$. Rows 4 and 5 are subcritical, with the hard onset of
Remark~\ref{rem:finite}. The table is a selection from a larger set. Of the
$4\,982$ four-reaction networks, $4\,298$ satisfy the first criterion with
$Q_1$ of the same strict sign, and $153$ of these satisfy the third. Among
those $153$, $44$ have $m_{01}=0$. All of these are supercritical, and $42$
are of activator--inhibitor type. Row 4 is the only subcritical one with
$\eta\ge1$.

\begin{table}[ht]
\centering
\footnotesize
\begin{tabular}{@{}rlccccc@{}}
\toprule
& network & mechanism & label & $\theta\to1$ & $\eta$ & Hopf up to\\
\midrule
1 & $0\to\ce X,\ \ce Y\to0,\ \ce X+\ce Y\to2\ce Y,\ 2\ce X\to3\ce X$
  & AI & super & super & $\infty$ & $\eps\ge1$\\
2 & $\ce Y\to\ce X+2\ce Y,\ \ce X\to0,\ \ce X+\ce Y\to2\ce X,\ 2\ce X\to3\ce X$
  & SD & super & super & $\infty$ & $\eps\ge1$\\
3 & $0\to\ce Y,\ \ce Y\to\ce X,\ \ce X\to0,\ 2\ce X+\ce Y\to3\ce X$
  & SD & super & super & $1$ & $\eps\ge0.2$\\
4 & $0\to2\ce X+\ce Y,\ \ce X\to0,\ \ce X+\ce Y\to2\ce X,\ 2\ce X\to3\ce X$
  & SD & sub & sub & $1$ & $\eps=0.68$\\
5 & $\ce Y\to\ce X,\ \ce X\to\ce X+\ce Y,\ \ce X+\ce Y\to0,\ 2\ce X\to3\ce X$
  & AI & sub & sub & $1/9$ & $\eps=0.104$\\
\bottomrule
\end{tabular}
\caption{Robust minimal canard networks (X fast, Y slow). Label: the
criticality for all rate constants; $\theta\to1$: the label of $Q_1$, the
limit of strong storage of X. $\eta=-n_{01}m_{10}/m_{01}^2$, the estimate
$\eps_{\mathrm{BT}}$ for a canard point with $x=y$. Hopf up to: the end of
the Hopf branch, continued from $\eps=10^{-5}$, with the canard point at
$(1,1)$ and $k=v$; ``$\ge$'' means that the branch was followed up to that
value (for row 3 the continuation failed between $0.2$ and $0.5$). Row 3 is
the autocatalator.}
\label{tab:robust}
\end{table}

\paragraph{Molecularity and thermodynamics.}
Corollary~\ref{cor:coop} shows that molecularity three is necessary rather
than an artefact of the enumeration: some step that produces X must have a
rate at least quadratic in X, and no network whose complexes have at most two
molecules has a canard point. This is why the classical planar models with
canard explosions, such as the autocatalator and the Brusselator, are cubic.
In some networks of the list, molecularity three occurs only among the
products. Counting X and Y only, $4\,974$ minimal networks have no reactant
with three molecules; in $658$ of them the only complex with three molecules
is the product of $2\ce X\to3\ce X$, and $27\,096$ further networks have
exactly one such reactant. Cubic steps are usually effective descriptions of
bimolecular steps through intermediates. By Corollary~\ref{cor:bimol},
resolving them through fast dimers keeps the canard points and changes the
criticality only through the term $2\theta n_{01}m_{10}n_{20}$. Every
nondegenerate canard point keeps its criticality under sufficiently weak
dimerisation, and the one-sided labels of $57\,796$ networks hold for every
strength. For other networks, the opposite criticality can appear under
arbitrarily weak dimerisation.

Canard explosions need a system driven away from equilibrium.
Complex-balanced mass-action systems, and in particular detailed-balanced
ones, have a strict Lyapunov function and hence no limit cycles
\citep{HornJackson1972}. The networks of the classification are open systems
in which constant pools are absorbed into the rate constants, as in the
Oregonator model of the Belousov--Zhabotinsky reaction. Irreversibility, on
the other hand, is not essential. If a minimal network with $\rank B=3$ is
completed by the reverse of each of its reactions, then by
Proposition~\ref{prop:inherit} the larger network has canard points of the
same criticality, at small reverse fluxes.

\paragraph{The Oregonator.}
The three-variable Oregonator \citep{FieldNoyes1974} has no conservation law.
With bromate and the organic substrate absorbed into the rate constants, its
reactions change HBrO$_2$, bromide and the oxidised catalyst independently.
So it is not planar on stoichiometric classes. With HBrO$_2$ and bromide both
fast, it has two fast species, a case that remains open (see
\emph{Extensions}). It does reach the present class through the two-variable
reduction of \citet{TysonFife1980}, in which bromide is eliminated as the
fastest species:
\[
  \eps\dot u=u-u^2-fv\,\frac{u-q}{u+q},\qquad \dot v=u-v ,
\]
where $u$ is HBrO$_2$ and $v$ the oxidised catalyst. Multiplication by
$u+q>0$ changes time but neither the orbits nor the canard points nor the
criticality. It gives a mass-action system of the form~\eqref{eq:sys} with
$\ce X=u$ fast and $\ce Y=v$ slow, namely the network
\begin{gather*}
  \ce X\to2\ce X,\ 2\ce X\to3\ce X,\ 3\ce X\to2\ce X,\ \ce X+\ce Y\to\ce Y,\
  \ce Y\to\ce X+\ce Y,\\
  \ce X\to\ce X+\ce Y,\ 2\ce X\to2\ce X+\ce Y,\ \ce X+\ce Y\to\ce X,\ \ce Y\to0 ,
\end{gather*}
with rate constants $q,1-q,1,f,fq$ for the fast reactions in the first line
and $q,1,1,q$ for the slow ones in the second, and molecularity three. So the two-variable Oregonator belongs to the class classified here.
Its network contains six minimal canard networks of the list: two
\emph{super}, one \emph{sub} and three with $A\equiv0$. Its canard points lie
on $u=v$ at the two folds of the fast nullcline. As $q\to0$ they tend to
$f=1+\sqrt2$ and $f=1/2$, the classical ends of the oscillatory range, and
both are of activator--inhibitor type. At the lower fold $f_{xx}>0$ and
$Q<0$, so the explosion is supercritical. At the upper fold $f_{xx}<0$, made
possible by $3\ce X\to2\ce X$ (Remark~\ref{rem:orient}), and $Q>0$, so it is
subcritical. The first Lyapunov coefficient confirms both signs for $q$ from
$1/40$ to $2\cdot10^{-4}$ and $\eps$ from $10^{-3}$ to $10^{-5}$. The
elimination of bromide is itself a singular reduction, and the criticality
of the three-variable model can differ from that of the reduction. Let
bromide relax on the time scale $\delta$, and take $q=8\cdot10^{-4}$ and
$\eps=10^{-3}$. Near $f=1/2$ the singular Hopf bifurcation is subcritical
for every $\delta$ from $10^{-8}$ to $10^{-2}$. Near $f=1+\sqrt2$ it is
supercritical for $\delta\le3\cdot10^{-6}$, as the reduction predicts, but
subcritical for $\delta\ge10^{-5}$. At that canard point $u\approx(1+\sqrt2)\,q$, and
bromide relaxes at the rate $(q+u)/\delta$. So, heuristically, the reduction
needs $\delta\ll\eps(q+u)\approx3\cdot10^{-6}$, which is far more demanding
than $\delta\ll\eps$.

\paragraph{Extensions.}
The classification concerns two species, but its networks are building
blocks for larger ones, and the effect of further species depends on their
speed. Weakly coupled species that relax on the fast time scale leave canard
points and their criticality intact (Remark~\ref{rem:species}). Species that
are faster than X make the reduction singular, and terms from the eliminated
species can enter $A$. For fast storage of X, as in the dimers of
Section~\ref{sec:bimol}, this term is explicit
(Proposition~\ref{prop:storage}); fast species that are not at
equilibrium with X remain to be treated. For networks with more species but one slow direction,
conservation laws make the dynamics planar on each stoichiometric class. In
forthcoming work we show that the flux criterion extends to this case. When
one species is fast, the cubic becomes a weighted sum of the planar cubics of
this paper, one for each pair formed by X and a slow species, with weights
that depend on the concentrations of the slow species, and the mode
criterion with its bound of three modes carries over. Three directions
remain open: several slow directions, where folded singularities organise
mixed-mode oscillations \citep{SzmolyanWechselberger2001,Desroches2012}; two
or more fast species, where the fold condition is no longer linear in the
fluxes; and splits in which some changes of Y are fast. The flux viewpoint
and the mode criterion should carry over to all three. In each case the sign
conditions become conditions on more linear forms, and Carathéodory's
theorem bounds the number of modes needed.

\paragraph{Data and code availability.}
The complete list of minimal canard networks with their labels, the
elementary modes, the certificates of the labels (including those in the
limit of strong dimerisation), the logs of all computations and checks, the
Julia code that reproduces the classification, and the independent verifier
in Python are openly available on Zenodo
\citep{CanardPlanarData2026}, at
\url{https://doi.org/10.5281/zenodo.23698722}. The README of the record
describes the file formats, gives examples of queries on the list, and
explains how to run the Julia pipeline and the independent verifier. The
computations for Table~\ref{tab:robust} and for the Oregonator are included.

\bibliographystyle{plainnat}
\bibliography{refs}

\end{document}